\documentclass[a4wide]{article}

\usepackage{graphicx}
\usepackage{multirow}
\usepackage{amsmath,amssymb,amsfonts}
\usepackage{amsthm}
\usepackage{mathrsfs}
\usepackage[title]{appendix}
\usepackage{xcolor}
\usepackage{textcomp}
\usepackage{manyfoot}
\usepackage{booktabs}
\usepackage{algorithm}
\usepackage{algorithmicx}
\usepackage{algpseudocode}
\usepackage{listings}
\usepackage{subcaption}

\theoremstyle{thm}

\begin{document}

\title{How different is rational approximation from piecewise polynomial approximation?}

\author{Daan Huybrechs}

\maketitle

\abstract{The first aim of this paper is to show that there is merit to the question posed in the title. Indeed, for certain function classes, approximation by rational functions and by piecewise polynomials are surprisingly similar. Between these two, polynomials are more widely studied and more widely used. One setting in which both approaches are equally well understood, and equally practical in use, is that of approximating univariate functions with point singularities. In this context we can fully address the question. We review classical literature on the topic which shows that both approaches do indeed achieve similar convergence rates. However, rational approximations come with significantly smaller constants. Owing to recent advances in practical rational approximation, we can augment the discussion with a comparison of modern numerical techniques that achieve the optimal rates. We end by showing numerically that the difference becomes even more pronounced in several variables.}

\section{Introduction}\label{s:intro}

Polynomials are ubiquitous in numerical analysis and scientific computing. Their role is usually one of approximation, as it is easier in algorithms to manipulate polynomials than it is to manipulate other types of continuous functions. Moreover, approximation by polynomials is well understood and this has been the case for more than a century. It was realized early on that the rate of convergence of polynomial approximation is related to the smoothness of the function involved: the smoother the function, the faster the convergence. This insight grew out of the understanding of polynomial approximations to the non-smooth function $f(x)=|x|$, by de la Vall{\'e}e Poussin and Bernstein in the early twentieth century, which played a seminal role in the development of approximation theory~\cite{deLaValleePoussin1908,bernstein1912}. Their choice was inspired by Lebesgue who, in his first published paper~\cite{Lebesgue1898}, showed that the Weierstrass approximation theorem can be reduced to approximating $|x|$. While this example served well to illuminate the connection between derivatives and convergence rates, it also highlighted from the very beginning one of the main deficiencies of polynomial approximation. The polynomial approximation of the non-smooth function $|x|$ itself is, in fact, rather poor.

The search for alternatives started quickly. If polynomials are indeed ubiquitous in scientific computing nowadays, it is possibly not because they approximate smooth functions well, but because piecewise polynomials approximate much broader classes of functions well. Splines and B{\'e}zier curves in Computer Aided Geometric Design (CAGD) are piecewise polynomials. Most basis functions employed in finite element methods (FEM) and boundary element methods (BEM) are piecewise polynomials. Neural networks employing the ReLU activation function represent piecewise polynomials. These three examples encompass a very broad range of applications. In many such settings it is very restrictive to demand of a function to be globally smooth. In contrast, most functions one encounters seem at least piecewise smooth. Relevant exceptions aside, such as stochastic quantities and noisy observables, for the purpose of approximation piecewise polynomials are a very sensible generalization of polynomials.

From algebraic and mathematical perspectives, rational functions are arguably a more natural generalization to consider. Rational functions are easily defined as quotients of polynomials, $r(x) = p(x) / q(x)$, and as such they inherit and extend a lot of their algebraic properties. Rational functions feature in basic calculus through partial fractions, often to compute integrals, and they shine in complex analysis, system theory and many other fields. However, in approximation theory, their role has been muted in comparison with that of piecewise polynomials. The most prominent counterexample are Pad{\'e} approximants, which Pad{\'e} explicitly devised for approximation in the late nineteenth century as a student of Hermite~\cite{Pade1892}. The study of Pad{\'e} approximants throughout the twentieth century revealed that, unlike polynomials, rational functions can approximate functions with poles and branch points~\cite{Brezinski1991,Bultheel1987}. The development and study of rational functions cannot be dismissed. However, it did not reach the same scale as that of piecewise polynomials.

It is difficult to quantify or compare the amount of research effort invested in a field. Foregoing any claim of rigor, we resort to some anecdotal observations. In 1998 a review paper on nonlinear approximation by DeVore appeared in Acta Numerica~\cite{DeVore1998}. In this 100 page paper, DeVore devotes precisely one page to rational functions. Intriguingly, \emph{\S 6.6 Rational approximation} is a subsection of \emph{\S 6 Piecewise polynomial approximation}. There can be no doubt on the merit of the title of the current paper, as DeVore explicitly writes: \emph{``The status of rational approximation is more or less the same as for piecewise polynomials.''} We recall the theoretical context on which this statement is based further on. It was made not for lack of enthusiasm for rational functions, on the contrary: \emph{``Thus, on the one hand the approximation problem is solved but on the other hand the news is somewhat depressing since there is nothing to gain or lose (in the context of the approximation classes) in choosing rational functions over piecewise polynomials.''} Expectations about rational functions had been different at the time, but analysis clearly showed similar convergence rates for similar function classes between rational functions and piecewise polynomials.

The theme of the current paper is illustrated well by another quote from the same review paper~\cite{DeVore1998}, still in \S6.6: \emph{``There have been several other important developments in rational approximation. One of these was Newman’s theorem (\ldots) which showed that the function $f(x) = |x|$ satisfies $r_n(f)_\infty = e^{-c \sqrt{n}}$ (a very stunning result at the time).''}\footnote{The quantity $r_n(f)$ here is Newman's notation for the maximum pointwise error on $[-1,1]$. We formally define the symbol further on in \S\ref{ss:equality}.} Indeed, in~\cite{Newman1964} Newman had proved exponential convergence of rational functions in the maximum norm on $[-1,1]$ to $|x|$, the very same function that de la Vall{\'e}e Poussin and Bernstein had chosen to base their exposition of polynomial approximation on. There is additional evidence of the sentiment of surprise. H. Stahl writes in~\cite{stahl1993}: \emph{``While Bernstein's investigations on best polynomial approximations of $|x|$ and $|x|^\alpha$ were published in the period between 1909 and 1938, the study of best rational approximation of $|x|^\alpha$ was only started in 1964 by Newman's surprising (at the time) result that\ldots''}.

Newman himself did not seem to make much of the importance of the absolute value function for rational approximation, writing: \emph{``While it is true that $|x|$ loses much of its special significance, it is nevertheless of some interest to determine the approximation to $|x|$ by $n$-th order rational functions.''} However, he did express his own surprise in the next sentence: \emph{``Now it is known that, in some overall sense, rational approximation is essentially no better than polynomial approximation, and this suggests the naive guess that the order of approximation of $|x|$ by $n$-th order rational functions is also $1/n$.''} If $1/n$ convergence was indeed the prevailing intuition at the time, then exponential convergence is certainly startling.

Rational approximations do much better than polynomial approximations~\cite{Trefethen2024Laplace}. However, they are not better - in the quantification used in~\cite{DeVore1998} and references therein - than piecewise polynomials. This leads to the question: how different are they?

\subsection{Scope and goal of the paper}

The current paper grew out of a talk presented in April 2025 in Banff, Canada, as part of the BIRS research workshop `Challenges, Opportunities and New Horizons in Rational Approximation'. It took over fifty years from Bernstein's work to Newman's paper to appreciate the value of rational approximation for non-smooth functions. In the decades that followed, several numerical algorithms for rational approximation were developed. The recent acceleration of these efforts led to the workshop, but also allows us to juxtapose modern approximation methods with earlier ones.

We cannot continue to use $|x|$ as an example, since that function is itself piecewise polynomial. It would be like comparing polynomial and rational approximation for the Runge function $1/(1+x^2)$. However, polynomial and rational approximation of $|x|$ on $[-1,1]$ is equivalent to approximation of $\sqrt{x}$ on $[0,1]$. Indeed, since $|x|$ is symmetric on $[-1,1]$, so are its best polynomial approximation $p$ and rational approximation $r$. Thus, $p$ and $r$ are functions of $x^2$, which means that $\tilde{p}(x) = p(\sqrt{x})$ and $\tilde{r}(x) = r(\sqrt{x})$ are also polynomial and rational functions of $x$, respectively. The maximum error on $[0,1]$ remains unchanged since, with $t=\sqrt{x}$,
\[
 \sqrt{x} - \tilde{p}(x) = \sqrt{x} - p(\sqrt{x}) = t - p(t), \quad x \in [0,1].
\]
Newman did not see $|x|$ as a very relevant special case. By quoting results for the more general setting of $x^\alpha$ with $\alpha > 0$, we will see in hindsight that the case $\alpha=1/2$ of the square root function actually does convey several general principles.

\subsection{Structure of the paper}

We briefly recall the theoretical setting of~\cite{DeVore1998} in \S\ref{s:approx}. For piecewise polynomials, we recall the theoretical results by DeVore and Scherer, and Babu{\v s}ka and collaborators, on so-called \emph{free knot splines} in \S\ref{s:piecewise}. We do the same for rational approximation in \S\ref{s:rational} with historical, and in this case also recent, references. All results are illustrated numerically in~\S\ref{s:numerical}, including an example with two and three variables. We end with a brief discussion in~\S\ref{s:discussion}.

\section{How similar is rational approximation to piecewise polynomial approximation?}\label{s:approx}

\subsection{Linear and nonlinear approximation}

We start with the observation that nodes and poles play similar roles.

By linear approximation schemes we mean all schemes in which $f(x)$ is approximated by a linear combination of basis functions $\{\phi_k(x)\}_{k=1}^n$, such that all coefficients $a_k$ in $f(x) \approx \sum_{k=1}^n a_k \phi_k(x)$ are linear functionals of $f$. Nonlinear approximation is everything else and includes in particular any form of adaptivity.

Approximation by piecewise polynomials is linear if the nodes or knots, the points of possible discontinuity, are prescribed. Approximation can be based on interpolation, least squares, projection or any other linear method. Methods vary in the amount of smoothness that is retained at the nodes. A lot of flexibility, especially in multivariate settings and complicated geometries, is achieved by not enforcing such smoothness at all. This underlies the success of discontinuous Galerkin methods~\cite{HesthavenWarburton2008}, although one should add that numerical flux conditions certainly render this setting more complicated than simple piecewise approximation. The approximation becomes nonlinear if the nodes are allowed to vary, sometimes referred to as free-knot or variable knot splines. It happens in CAGD when knot locations are used as design parameters (see, e.g., \cite{PieglTiller1997}), in FEM when a mesh is adaptively refined and in neural networks when weights and biases are updated. In all these nonlinear settings, finding the `best' approximation can be challenging. There is no standard way to represent piecewise polynomials either with or without prescribed nodes.

Approximation by rational functions is linear if the poles are prescribed. Numerically, in this setting rational functions are best utilised in partial fraction form,
\begin{equation}\label{eq:partial_fraction}
 f(x) \approx \sum_{k=1}^n \frac{a_k}{x-p_k},
\end{equation}
possibly with a normalization factor if a pole $p_k$ is close to the approximation domain. Such linear rational approximation is the basis of the recent lightning-class of methods for the solution of elliptic PDEs in two dimensions~\cite{GopalTrefethen2019b,GopalTrefethenPNAS2019}. Approximation by rational functions is nonlinear if the poles are allowed to vary. In this setting, an emerging representation format in algorithms such as Adaptive Antoulas Anderson (AAA)~\cite{nakatsukasa2018aaa} is the barycentric form~\cite{berrut2004barycentric},
\begin{equation}\label{eq:barycentric}
 f(x) \approx \frac{\sum_{j=1}^n \frac{\alpha_j}{x-x_j}}{\sum_{j=1}^n \frac{\beta_j}{x-x_j}}, \quad \mbox{or} \quad f(x) \approx \frac{\sum_{j=1}^n \frac{w_j f(x_j)}{x-x_j}}{\sum_{j=1}^n \frac{w_j}{x-x_j}}.
\end{equation}
The second form achieves interpolation of $f$ in the \emph{support points} $x_j$ (not to be confused with poles of the approximation). Here, finding suitable \emph{barycentric weights} $w_j$ is a nonlinear problem, but the formulation is linear in the function values $f(x_j)$ and thus in $f$. Indeed, the weights determine the poles as the roots of the denominator, which is independent of $f$.

\subsection{Localisation and approximation order}

Piecewise polynomials have an intrinsically localized nature, especially if no smoothness across nodes or mesh boundaries is enforced. In that case the polynomials on each piece are entirely independent. It is clear that, similar to wavelets, functions can be well approximated by piecewise polynomials if they are locally smooth. Associated with local polynomials is a local notion of order, related to the degree of the polynomial. Piecewise polynomials are free to vary order locally, depending on local smoothness of the function. In practice, there is often a trade-off between reducing approximation error by refining the pieces or by increasing the local polynomial degree, or both. This trade-off can be challenging to resolve.

A rational function is meromorphic and hence cannot have compact support. Moreover, there is no local notion of order, only a global rational degree. It may therefor not be immediately obvious that rational functions can successfully approximate functions with locally varying smoothness. However, the mechanism at play is exactly the same and is revealed through several perspectives. We elaborate in~\S\ref{ss:equality} below on the theoretical analysis of smoothness spaces: rational functions behave (or rather, can behave) exactly like adaptive, locally refined polynomial approximations. In~\S\ref{ss:resolution} we offer the perspective of the resolution of point singularities, which we explore in depth in later sections of the paper. One could say that a rational function can think globally but act locally.

Reasoning in the other direction, from local to global, a piecewise polynomial has no global properties. Their approximation properties do not extend beyond the approximation domain itself and they have no meaningful extension into the complex plane. Rational functions have this in the extreme: Pad{\'e} approximants matching Taylor series to all orders at a single point can reveal distant singularities of a function.

\subsection{Equality of approximation spaces}\label{ss:equality}

We recall the main technical result as formulated in~\cite{DeVore1998} and attributed to Pekarski~\cite{Pekarski1986} and to Petrushev~\cite{Petrushev1988}. The precise statement is~\cite[(6.30)]{DeVore1998}:
\begin{equation}\label{eq:equality}
    {\mathcal A}_q^\alpha(L_p(\Omega),(\mathcal{R}_n)) = {\mathcal A}_q^\alpha(L_p(\Omega),(\Sigma_{n,r})), \quad 0 < \alpha < r.
\end{equation}
Here, $\Omega$ is a bounded interval on the real line. Both the left and right hand sides are sets of functions, subsets of the standard Lebesgue space $L_p(\Omega)$, associated with:
\begin{itemize}
    \item $\mathcal{R}_n$ (left): univariate rational functions up to degree $n$, meaning that numerator and denominator are polynomials of degree at most $n$,
    \item $\Sigma_{n,r}$ (right): piecewise polynomials of fixed degree smaller than $r$ with $n$ free pieces.
\end{itemize}
Approximation by rational functions is measured by the error~\cite[(6.29)]{DeVore1998}
\[
r_n(f)_p := \inf_{R \in {\mathcal R}_n} \Vert f - R \Vert_{L_p(\Omega)}.
\]
Similarly, the approximation error associated with piecewise polynomials is~\cite[(6.18)]{DeVore1998}
\[
 \sigma_{n,r}(f)_p := \inf_{S \in \Sigma_{n,r}} \Vert f - S \Vert_{L_p(\Omega)}.
\]
The approximation space in the left hand side of~\eqref{eq:equality} contains all functions that can be approximated at a rate $n^{-\alpha}$ by rational functions, in the sense that
\[
{\mathcal A}_q^\alpha(L_p(\Omega), ({\mathcal R})_n) = \left\{ f \in L_p(\Omega) : \left( \sum_{n=1}^\infty (n^\alpha r_n(f)_p)^q \frac{1}{n} \right)^{1/q} < \infty\right\},
\]
with a standard modification for the case $q=\infty$ (the rate remains $n^{-\alpha}$ regardless of $q$). There is a similar definition for the right hand side in terms of $\sigma_{r,n}$, which comes with the restriction $0 < \alpha < r$.

Within these strict definitions, there is no difference in approximation rates between piecewise polynomials and rational functions. That is precisely what equality~\eqref{eq:equality} expresses: the sets of functions that give rise to a certain convergence rate are exactly the same in both settings. This is the basis for the claim that the status of rational approximation is close to that of piecewise polynomials. It is provably correct, although of course subject to the conditions of the statement.

\subsection{Which one is more efficient?}\label{ss:efficient}

The method of proof in the given references above~\cite{Pekarski1986,Petrushev1988}, as well as in earlier and later work such as a multivariate extension in~\cite{DeVoreYu1986}, typically starts from a piecewise polynomial approximation. A rational function is constructed explicitly by a weighted combination of the polynomial pieces, multiplied by rational cut-off functions. Since statement~\eqref{eq:equality} is only asymptotic, it suffices to show that a rational function exists with a degree that is linear in the number of pieces, with cut-offs chosen not to impact the algebraic convergence rate.

In principle, such proofs are constructive. However, they would lead one to think that piecewise polynomials are more efficient. The rational functions produced in these constructions are more complicated than the piecewise polynomial approximation and, hence, the polynomial approach seemingly achieves the same convergence rate with fewer degrees of freedom. We shall show in the numerical experiments later on that the situation is, in fact, exactly the opposite. For the same convergence rate, rational functions require (much) fewer degrees of freedom.

Such a comparison does depend on how polynomials are represented. For piecewise polynomials of fixed degree, which is the setting of~\eqref{eq:equality}, the use of basis functions such as B-splines requires fewer degrees of freedom than the independent representation of the polynomials. The construction of efficient representations for the setting of varying order is much more involved, as are splines in multivariate settings.

\subsection{Resolution of singularities by varying order}\label{ss:resolution}

Any comparison should not end with the fixed-degree setting of~\eqref{eq:equality}. Piecewise polynomials with varying degree have considerably improved approximation properties and have been the subject of intense research. The research effort originated in the pursuit of hp-adaptivity, with precise results being developed early on for functions with point singularities (reviewed, e.g., in~\cite{BabuskaSuri1994}).

Here, too, the comparison to rational functions remains apt. In fact, as illustrated with several more examples in~\cite{huybrechs2024sigmoids}, rapidly converging approximation schemes for the resolution of singularities seem to employ universal mechanisms. For piecewise polynomials with varying order and likewise for rational functions, convergence rates have been studied in detail. Both converge at \emph{root-exponential rates}, i.e., rates of the form $e^{-c \sqrt{n}}$. However, in this setting, their exact rates differ. We recall the known results in the next two sections.

\section{Piecewise polynomial approximation to $x^\alpha$ on $[0,1]$}\label{s:piecewise}

The approximation of $x^\alpha$ on exponentially graded meshes was studied by DeVore and Scherer in~\cite{DeVoreScherer1980,Scherer1981} and, in the context of the solution of ODEs and PDEs, by Babu\v{s}ka and collaborators~\cite{GuiBabuvska1986a,GuoBabuvska1986}. We consider only bounded functions, i.e. $\alpha > 0$.

For historical reasons we focus on approximations with respect to the maximum norm. That is, for an $n$-dimensional set of functions $V_n$ we aim to minimize
\[
E_n(f) = \min_{v \in V_n} \max_{x \in [0,1]} | f - v|.
\]
For a more complete analysis involving other norms we refer to~\cite{Schwab1998,Melenk2002}.

Optimal approximations are found by balancing the local errors. There is an approximation error on each interval of the graded mesh. In addition, there is a truncation error close to the singularity.

\subsection{Root-exponential convergence by exponential clustering}

Exponentially refined meshes on the interval $[0,1]$ can be formulated in terms of a parameter $0 < \eta < 1$ with nodes $t_k = \eta^k$, for $k = 0,\ldots,m-1$, augmented with $t_m = 0$. Polynomials of degree $n_k$ are used on each interval $[t_{k+1},t_k]$. It is customary to approximate the function by $0$ on the interval $[0,t_{m-1}] = [0,\eta^{m-1}]$. This implies that the maximum error on that interval is given by $t_{m-1}^\alpha = \eta^{\alpha(m-1)}$. If all errors are balanced, this truncation error captures the convergence rate.

Standard polynomial approximation theory for analytic function quickly shows that exponentially graded meshes have uniform exponential convergence rates. For approximation on $[-1,1]$ the rate is $\rho^n$, where $\rho$ is the radius of the largest Bernstein ellipse in which $f$ is analytic~\cite{trefethen2019atap}. The interval $[\eta^{k+1},\eta^k]$ is mapped to $[-1,1]$ by
\[
m(x) = 2\frac{x - \eta^{k+1}}{\eta^k - \eta^{k+1}} - 1.
\]
For each interval, the nearest singularity is the branch point at $x=0$. That singularity is mapped to
\[
 m(0) = 2\frac{- \eta^{k+1}}{\eta^k - \eta^{k+1}} - 1 = -2 \frac{\eta}{1-\eta} - 1 = \frac{\eta + 1}{\eta - 1},
\]
which is indeed independent of $k$. The exponential rate of convergence on each piece is
\[
 \rho = |m(0)| + \sqrt{m(0)^2-1} = \frac{1+\sqrt{\eta}}{1-\sqrt{\eta}} > 1.
\]

These observations are sufficient to demonstrate root-exponential convergence in $n$. We can choose fixed degrees $n_k = m$, leading to $n = m^2$ degrees of freedom in total. Increasing $m$ by $1$ adds an interval close to $0$, while leaving the other intervals unchanged, and decreases the error everywhere by a constant factor $\rho$. Thus, we have exponential convergence in $m$ and, since $m \sim \sqrt{n}$, root-exponential convergence in $n$.

We can still optimize the choice of $\eta$. DeVore and Scherer conclude in~\cite{DeVoreScherer1980} that $\eta = (\sqrt{2}-1)^2 \approx 0.172$ is an optimal choice for all $0 < \alpha < 1$.\footnote{This value was also obtained by Gonchar earlier on~\cite{Gonchar1972Piecewise}, but his overall convergence rate was not optimal.} If we write $\eta = \delta^2$, then the choice $\delta = \sqrt{2}-1$ happens to be such that
\[
 \rho = \frac{1+\sqrt{\eta}}{1-\sqrt{\eta}} = \frac{1+\delta}{1-\delta} = \delta.
\]
In the proof, $\delta$ arises as the value satisfying this property and achieving the balance.

\subsection{Tapering with a linear degree vector}

We have yet to describe the optimal degrees. As it turns out, they are linearly decreasing within subintervals approaching $x=0$. Degrees closer to the singularity can be smaller than degrees further away. Intuitively, this is simply because $x^\alpha$ is smaller near $x=0$, while the subintervals also become smaller: the absolute error criterion loosens by a fixed factor with increasing $k$~\cite[\S3]{huybrechs2024sigmoids}. A consequence of this is that, while the grid points are exponentially distributed, the distribution of the degrees of freedom is slightly skewed, as there are fewer of them near the singularity. A similar effect is called \emph{tapering} in~\cite{TrefethenNakatsukasaWeideman2021} in the context of rational approximation.

The precise choice of degrees determined in~\cite{DeVoreScherer1980} is:
\begin{equation}\label{eq:degrees}
 n_k = \left\{ \begin{array}{cc} 1, & \mbox{if~} \alpha (m-k) < 1, \\ 
 \lfloor 2\alpha (m-k) \rfloor, & \mbox{otherwise.} \end{array}\right.
\end{equation}
Since the degree vector has linearly increasing elements, its sum is quadratic in $m$. More specifically, we have
\[
n \leq 1 + \alpha m (m+1).
\]
Hence, $m \sim \sqrt{n/\alpha}$, which means the size of the truncation error is
\begin{equation}\label{eq:optimal_rate_poly}
 t_{m-1}^\alpha = \eta^{\alpha (m-1)} \sim \eta^{\sqrt{\alpha n}} = \delta^{2 \sqrt{\alpha n}}.
\end{equation}
The approximation error has the same size. In an error estimate, the decaying factor $\rho^{-n_k}$ multiplies the size of the function on the Bernstein ellips~\cite[Theorem 8.2]{trefethen2019atap}. On the interval $[\eta^{k+1},\eta^k]$ that factor scales with $(\eta^k)^\alpha$. Thus, the approximation error is
\[
 \rho^{-n_k} \eta^{\alpha k} \sim \delta^{2\alpha (m-k)} \delta^{2 \alpha k} = \delta^{2 \sqrt{\alpha n}}.
\]

We have only sketched heuristic arguments. DeVore and Scherer derive the results in a more systematic manner and prove that the rate~\eqref{eq:optimal_rate_poly} cannot be improved (up to algebraic pre-factors in $n$). Note that the argument relies on the polynomials being independent on each of the pieces, since we take $n = \sum_{k=0}^{m-1} n_k $.




\section{Rational approximation to $x^\alpha$ on $[0,1]$}\label{s:rational}

Newman's construction for $|x|$ and $\sqrt{x}$ was generalized to $x^\alpha$ and thoroughly investigated. The sharpest results are formulated by Stahl in~\cite{stahl1993,Stahl1994} (citing earlier results by Ganelius~\cite{Ganelius1979} and Vyacheslavov~\cite{Vyacheslavov1981}). Stahl shows that for rational approximations of degree $n$, the minimax error satisfies
\begin{equation}\label{eq:rate_rational}
 \lim_{n\to\infty} e^{2\pi\sqrt{\alpha n}}E_n(x^\alpha) = 4^{1+\alpha} |\sin(\pi \alpha)|.
\end{equation}
This implies a convergence rate of $e^{-2\pi \sqrt{\alpha n}}$. Comparing to the rate $e^{2\ln(\delta) \sqrt{\alpha n}}$ of~\eqref{eq:optimal_rate_poly}, rational approximations are more efficient by a wide margin in this setting since
\begin{equation}\label{eq:comparison}
 2\pi \approx 6.28 > -2 \ln \delta \approx 1.75.   
\end{equation}

Yet, many strong similarities remain. In the analysis of Stahl, as well as in the perspectives offered in the other references below, root-exponential convergence arises as a balance between truncation error and discretization error. Moreover, the rational approximations exhibit strong localisation, with poles playing a similar role to nodes.

\subsection{Tapered distribution of poles}

Using potential theory arguments, Stahl describes the distribution of poles in the form of a density function~\cite[Theorem 2.2]{Stahl1994}. However, this result is not fully explicit in a small neighbourhood of the origin. Trefethen, Nakatsukasa and Weideman focus on that region and explicitly describe a tapered exponential distribution of poles in~\cite{TrefethenNakatsukasaWeideman2021}. Herremans extends those results to achieve the optimal convergence rate $e^{-2\pi \sqrt{\alpha n}}$ in~\cite{herremans2021masterthesis} with an explicit numerical scheme, further analyzed in~\cite{HerremansHuybrechsTrefethen2023}.

Near-optimal rational approximations exist with a small number of large poles or poles at infinity that represent a smooth part of the function. The number of such poles is proportional to $\sqrt{n}$. In addition, a singular part can be resolved by finite poles on the negative real line given by the explicit formula
\begin{equation}\label{eq:poles}
 p_j = -C e^{-\sigma (\sqrt{n}-\sqrt{j})}, \qquad 1 \leq j \leq n.
\end{equation}
Formulas of this form were pioneered in~\cite{TrefethenNakatsukasaWeideman2021} and capture the tapering effect with the right scaling. They describe points that cluster exponentially towards the origin in such a way that $\log p_j$ is proportional to $\sqrt{j}$. This is faster than a purely exponential grading of the form $-C \eta^j$, in which case $\log p_j$ is proportional to $j$. We say `faster' in the sense that those points of~\eqref{eq:poles} that are close to the origin are spaced further apart in comparison to geometric grading. Away from the origin the points look similar.

The optimal choice of $\sigma$ is $2 \sqrt{2} \pi$ for $\alpha=1/2$. More generally, it is~\cite[(3.2)]{HerremansHuybrechsTrefethen2023}
\begin{equation}
    \sigma = \frac{2\pi}{\sqrt{\alpha}}.
\end{equation}
That makes the smallest pole $p_n$ of~\eqref{eq:poles} on the order of $e^{-2\pi\sqrt{n/\alpha}}$. On that scale, the function $x^\alpha$ has the size $e^{-2\pi\sqrt{\alpha n}}$. This is also the optimal convergence rate. At first sight, there seems no explicit truncation at play here. However, a truncation error related to the smallest pole does appear in each available path of analysis: in the potential theoretical analysis of~\cite{stahl1993,Stahl1994}, in the Hermite integral formula in~\cite{TrefethenNakatsukasaWeideman2021} and in an integral representation in~\cite{HerremansHuybrechsTrefethen2023}. A notion of discretization error stems from the numerical approximation of the latter two types of integrals.

\subsection{Localisation of rational approximations}

A rational function is always global and it may not be clear how a pole can have a localized effect. An elegant and highly informative visualization of their localization properties in the resolution of singularities is given in~\cite{huybrechs2024sigmoids}.

One can infer the independent properties of poles of a rational function from its partial fractions form. In the case at hand, that form would be
\[
 r_n(x) = \sum_{j=1}^n a_j \frac{p_j}{x - p_j} + \sum_{j=1}^{d \sqrt{n}} b_j x^j.
\]
Each pole can be associated with a basis function $r_j(x) = -\frac{p_j}{x - p_j}$. In this sense, each pole `contributes' independently to the rational function. Here, we have chosen to normalize such that $r_j(0)=1$. We note in addition that $r_j(x)$ tends to $0$ for large $x$.

The main theme of~\cite{huybrechs2024sigmoids} is the observation that different schemes for the resolution of a singularity on $[0,1]$ behave similarly when viewed after a logarithmic change of variables $s = \log x$. In the $s$ variable, the functions $r_j(x)$ have the shape of a sigmoid function. The sigmoid shape captures the switch of $r_j(x)$ from $1$ to $0$ and, still in the $s$ variable, it does so in a localized manner.

\begin{figure}[t]
  \centering
    \includegraphics[width=0.55\linewidth]{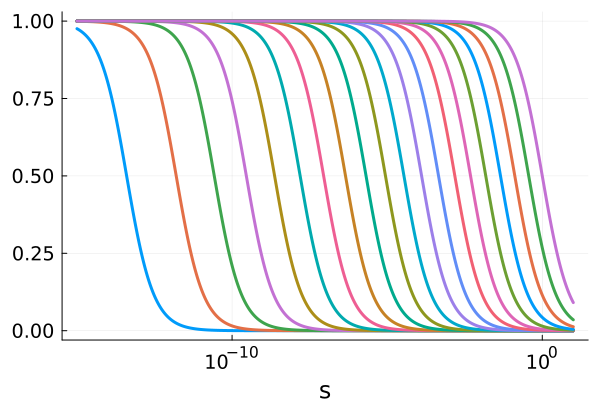}
    \caption{Illustration of the partial fraction basis functions $r_j(x) = - \frac{p_j}{x-p_j}$ as a function of $s = \log x$. The functions have a sigmoid shape which results in a clear localisation effect. In the transformed variable $s$, the approximation problem becomes a fit of a smooth function using translates of a smooth kernel.}
    \label{fig:sigmoid}
\end{figure}

The basis functions are shown in Fig.~\ref{fig:sigmoid}. Although the poles live on the negative real axis, their influence manifests itself near the intervals $[-p_{j+1},-p_j]$. This is strikingly similar to the intervals of a graded mesh. The localization is approximate, not highly accurate, and certainly not accurate enough to be seen as cut-off functions. Nevertheless, there is clearly a strong localisation effect at play.

\subsection{Efficiency and numerical stability}

The implementation and manipulation of piecewise polynomials is straightforward, at least when the polynomials of the pieces are independent.\footnote{Of course there is a rich body of literature on the topic of splines, but in multivariate settings or with variable order their implementation becomes very involved.} However, rational functions are catching up. If the perceived complexity of rational functions has played a role in the past, this should no longer be the case. The approximations we have discussed in this paper can be readily computed in a variety of ways. For this paper we have used approximations in partial fraction form~\eqref{eq:partial_fraction} by least squares fitting as in~\cite{GopalTrefethen2019b,GopalTrefethenPNAS2019,HerremansHuybrechsTrefethen2023}, and in barycentric form~\eqref{eq:barycentric} by AAA~\cite{nakatsukasa2018aaa}. Other recent developments in nonlinear rational approximation include rational Krylov fitting and related methods \cite{BerljafaGuttel2015,BerljafaGuttel2017}, as well as vector fitting, originally introduced in \cite{GustavsenSemlyen1999} and subsequently refined in \cite{DrmacGugercinBeattie2015Matrix,DrmacGugercinBeattie2015Quadrature}. Closely related approaches arising in model reduction and rational approximation include projected nonlinear least-squares methods \cite{HokansonMagruder2020} and more recent work on $L^2$ and $L^\infty$ rational approximation \cite{AckermannReiterTrefethen2025}.

Although least squares fitting does lead to ill-conditioned linear systems, their conditioning does not necessarily prevent high accuracy. Separately from the developments in rational approximation, insights in these phenomena grew out of the analysis of generalizations of bases to frames~\cite{adcock2019frames} and other redundant sets~\cite{herremans2026sampling}, settings in which ill-conditioning can arise as a relatively benign side-effect of the fact that a single function can have multiple representations.

\subsection{Optimality and non-optimality of the approximations}

The aim of the paper is to compare piecewise polynomials to rational functions, more so than to describe optimal approximations of $x^\alpha$. Still, since that is our main example, we include some more relevant results.

The complexity of approximating analytic functions on $(-1,1)$ with Lipschitz order $r$ at the endpoints has been studied by Burchard and H{\"o}llig in~\cite{BurchardHollig1985}. Their research is based on the study of n-widths, a concept that arises more generally when studying the best possible approximation to a class of functions from an $n$-dimensional linear space, optimizing over all such spaces (see~\cite[\S9]{DeVore1998} and \cite{Pinkus1985}). Burchard and H{\"o}llig show that approximations in linear spaces can converge at a rate no better than $e^{-\pi/2 \sqrt{r n}}$. This rate is achieved by Whittaker's cardinal series as demonstrated by Stenger~\cite{Stenger1976}. We can choose $r = \alpha$. However, we should be careful to compare to these results, as our setting of $x^\alpha$ on $[0,1]$ is different. First, we have only one singular point, whereas the study in~\cite{BurchardHollig1985} allows two. This means that one should at least double the values of $n$ in this paper. Second, and more importantly, the optimality statement does not cover nonlinear approximations. Unfortunately, recent research on nonlinear n-widths, involving manifolds rather than linear spaces, do not explicitly cover the setting of point singularities~\cite{DeVoreHowardMicchelli1989,CohenDeVorePetrovaWojtaszczyk2022}.

These differences aside, piecewise polynomials do not achieve the optimal rate, not even when optimized for $x^\alpha$. However, since the graded mesh is independent of $\alpha$, the space of piecewise polynomials is actually a linear space.

Intriguingly, rational approximations do much better than the stated optimal result, even when taking the doubling of $n$ into account. This is not a contradiction, but rather an illustration of the assumption of linearity. The rational approximation is specific to $x^\alpha$. In the context of analytic functions with point singularities away from the approximation domain, rather than on its boundary, the difference between linear and nonlinear (adaptive) rational approximation can be quantified by the famous $\rho$ vs. $\rho^2$ rate of convergence~\cite{Rakhmanov2016rho2}. In this setting adaptive approximations converge at twice the rate of linear approximations.\footnote{The constant in the convergence rate $e^{-2\pi \sqrt{\alpha n}}$ of~\eqref{eq:rate_rational} differs by a factor of $4$ with that of the optimal rate $e^{-\pi/2 \sqrt{rn}}$ of~\cite{BurchardHollig1985} when $r=\alpha$. A factor of $\sqrt{2}$ is explained by the doubling of $n$. A $\rho$ vs. $\rho^2$ effect would account for another factor of $2$ if applicable, but still leaves a difference of $\sqrt{2}$ unaccounted for.}

In conclusion, in our running comparison it is important to emphasize that piecewise polynomial approximations on $[0,1]$ retain their convergence rate for a broader set of functions than rational approximations do, unless the poles are kept fixed.

\begin{figure}[ht]
  \centering
    \includegraphics[width=0.55\linewidth]{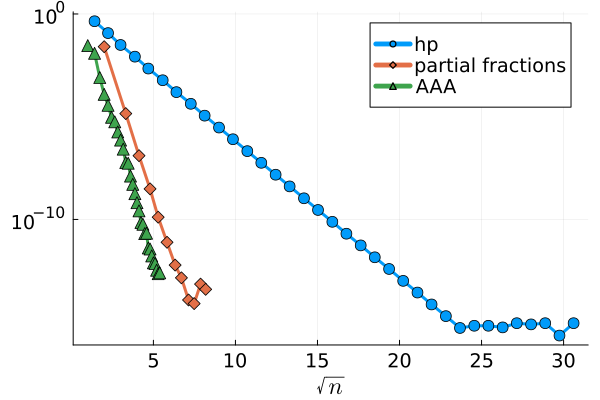}
    \caption{Illustration of approximation errors for $f(x)=\sqrt{x}$. The maximum error on $[0,1]$ is shown in logarithmic scale as a function of $\sqrt{n}$. A straight line implies root-exponential convergence. Results are shown until stagnation of the error.}
    \label{fig:logsqrt}
\end{figure}

\begin{figure}[ht]
  \centering
    \includegraphics[width=0.55\linewidth]{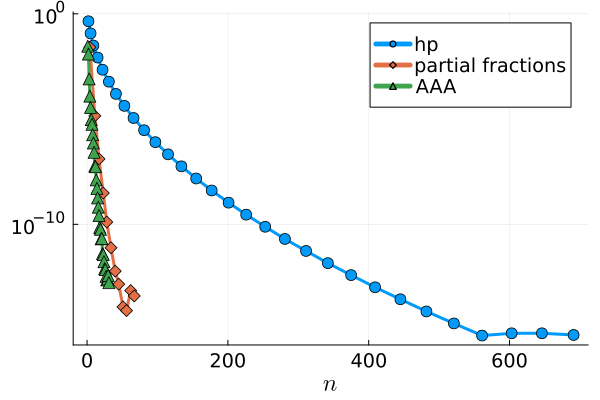}
    \caption{Similar to Fig.~\ref{fig:logsqrt} but for the function $f(x)=x^{3/4}$. The error is shown as a function of $n$ and displays root-exponential convergence.}
    \label{fig:loglinear}
\end{figure}

\section{Numerical experiments}\label{s:numerical}

The main experiment of the paper is a comparison between approximation of $x^\alpha$ on $[0,1]$ by piecewise polynomials and by rational functions.

\subsection{Illustration of optimal convergence rates in 1D}

We compare three schemes:
\begin{enumerate}
    \item Piecewise polynomial approximations with optimal choices of~\cite{DeVoreScherer1980} (recalled in \S\ref{s:piecewise}).
    \item Rational approximation in partial fraction form by least squares fitting, following the implementation described in~\cite{HerremansHuybrechsTrefethen2023}. We use~\eqref{eq:poles} with $\sigma = 2\pi/\sqrt{\alpha}$, $C=1$ and an additional polynomial part of degree $1.3\sqrt{n}$.
    \item The continuum AAA algorithm of~\cite{Driscoll2024continuumAAA} as implemented in the Julia package \texttt{RationalFunctionApproximation.jl}~\cite{RationalFunctionApproximation_jl}
\end{enumerate}
Results are shown in Fig.~\ref{fig:logsqrt} for $\alpha = 1/2$. All schemes achieve the convergence rates described in this paper. As expected, the convergence of the rational approximations is substantially faster than that of piecewise polynomials. Both rational approximations converge at the same rate, with continuum AAA having the better constant.

Results are also shown in Fig.~\ref{fig:loglinear} for $\alpha=3/4$. Here, we show the results as a function of $n$, rather than $\sqrt{n}$, in order to better appreciate the difference in the total number of degrees of freedom. Comparing around the $1e-10$ level of the error, piecewise polynomial approximations require approximately $8$ times more degrees of freedom ($226$ for hp versus $28$ for partial fractions).

\subsection{Approximation of a singular function in 2D}\label{ss:2d}

Having illustrated approximation behaviour in the univariate setting, we can present an outlook of what to expect in bivariate settings. The simplest experiment is an approximation via tensor-products of univariate approximants. Bivariate rational approximations for functions with curves of singularities have been described in~\cite{boulle2024rational}.

\begin{figure}[ht]
  \centering
    \includegraphics[width=0.55\linewidth]{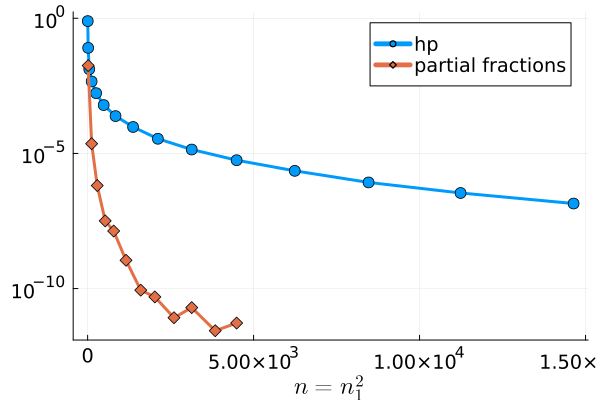}
    \caption{Maximum absolute error on $[0,1]^2$ in the approximation of $f_2$ as a function of the total number of degrees of freedom.}
    \label{fig:loglin_2d}
\end{figure}

\begin{figure}[ht]
  \centering
    \includegraphics[width=0.55\linewidth]{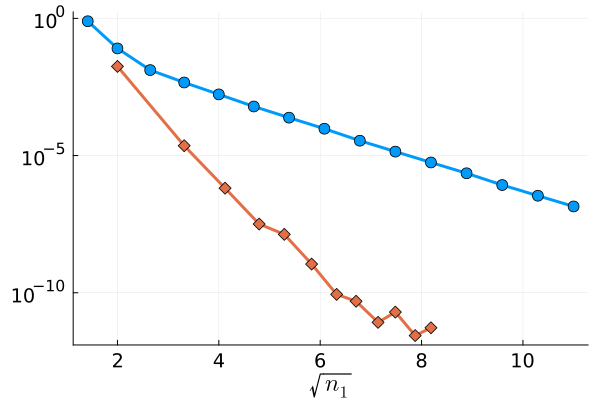}
    \caption{Similar to Fig.~\ref{fig:loglin_2d} but with the error shown in logarithmic scale as a function of $\sqrt{n_1}$. The straight lines confirm root-exponential convergence in this parameter.}
    \label{fig:logsqrt_2d}
\end{figure}

In order to show robustness of piecewise polynomial approximations and of rational functions, we choose to approximate the function
\[
f_2(x,y) = x^{y + \frac12} y^{x + \frac12}
\]
on the unit square $[0,1]^2$. This function has algebraic singularities along the edges $x=0$ and $y=0$, with the order of the singularity varying smoothly along the edge.

We reuse the setting of the univariate problem for $\alpha=1/2$ and extend everything via tensor products. Thus, we have $n = n_1^2$ degrees of freedom. For the rational approximation, we use Algorithm 2.1 of~\cite{boulle2024rational} with regularization tolerance set to $1e-15$.\footnote{Algorithm 2.1 presents an efficient solver for a least squares problem with a system matrix $A \otimes B$ having Kronecker structure. The solver constructs a regularized SVD solution based on the SVD of the factors $A$ and $B$ separately. Regularization is important in order to find solutions with modest norm.} We use the same algorithm without regularization for the piecewise olynomials. We expect to see convergence that is root-exponential in $n_1$ at a rate the corresponds to the choice $\alpha=1/2$. This is the order of the weakest algebraic singularity in the function $f_2$, achieved in the origin.

Fig.~\ref{fig:loglin_2d} shows the maximum absolute error on the unit square for piecewise polynomial and rational approximations. It is evident that the difference which existed in the univariate case has compounded: rational functions are more efficient by a wide margin. In Fig.~\ref{fig:logsqrt_2d} we plot the same errors as a function of $\sqrt{n_1}$. This plot confirms root-exponential convergence of both approximations in $n_1$, at rates that are similar to the optimal rate in one dimension for $\alpha=1/2$.

\subsection{Approximation of a singular function in 3D}\label{ss:3d}

We repeat the previous experiment, this time for the function of three variables
\[
f_3(x,y,z) = x^{y+z+1/2} y^{x+z+1/2} z^{x+y+1/2}
\]
on the unit cube $[0,1]^3$. As far as spectral approximation problems go, this function is quite challenging since it is singular in each point of the three faces $x=0$, $y=0$ and $z=0$ of the cube. The order of the algebraic singularity varies continuously.

For the computation of the approximations, we use a three-dimensional analogue of~\cite[Algorithm 2.1]{boulle2024rational}. For the rational approximation, we continue to use the regularisation threshold $1e-15$.

\begin{figure}[ht]
  \centering
    \includegraphics[width=0.55\linewidth]{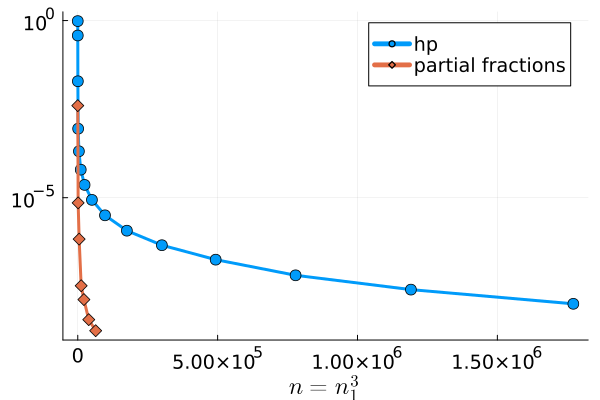}
    \caption{Maximum absolute error on $[0,1]^3$ in the approximation of $f_3$ as a function of the total number of degrees of freedom. The difference in convergence rates is stark.}
    \label{fig:loglog_3d}
\end{figure}

\begin{figure}[ht]
  \centering
    \includegraphics[width=0.55\linewidth]{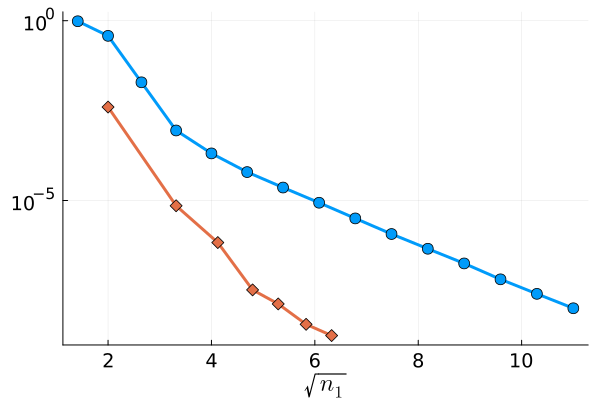}
    \caption{Similar to Fig.~\ref{fig:loglog_3d} but with the error shown in logarithmic scale on the Y-axis, as a function of $\sqrt{n_1}$ on the X-axis. The straight lines confirm root-exponential convergence in this parameter.}
    \label{fig:logsqrt_3d}
\end{figure}

The results are shown in Fig.~\ref{fig:loglog_3d} and Fig.~\ref{fig:logsqrt_3d}. Similar conclusions hold as in the previous experiment, with an even more pronounced difference between constants compared to the two-dimensional case. To reach an accuracy of about $1e-9$, the piecewise polynomials require $121^3=1,771,561$ degrees of freedom, while the rational function requires $28^3=21,952$ - a difference by a factor of around $80$.

\section{Concluding remarks}\label{s:discussion}

We have addressed the similarities between piecewise polynomials and rational functions mainly in the limited setting of approximating functions with an algebraic point singularity. Although some aspects of the comparison have indeed been specific to that setting, conceptually the comparison is certainly more general.

The study of rational functions of several variables is ongoing and rapidly developing. The analogy (and equality in certain contexts) between multivariate piecewise polynomials and rational functions has been studied by DeVore and Yu in~\cite{DeVoreYu1986}. Here, too, a multivariate piecewise polynomial can be approximated by a rational function. As in the univariate case this is proven using rational cut-off functions. Much more recently, effective algorithms for multivariate rational approximations have started to appear. The AAA algorithm using the barycentric representation was extended to multiple dimensions in the form of the p-AAA algorithm~\cite{CarracedoRodriguezBalickiGugercin2023,BalickiGugercin2026}. Multivariate approximations in partial fractions form for functions with curves of singularities were described in~\cite{boulle2024rational} and used to produce the numerical results in~\S\ref{ss:2d}.

The author wishes to thank the organizers of the BIRS workshop for bringing together a group of researchers with joined interests in rational functions and their many applications.

\bibliographystyle{abbrv}
\bibliography{bibliography}

\end{document}